\documentclass{article}

\newcommand{\N}{\ensuremath{ \mathbf N }}
\newcommand{\R}{\ensuremath{ \mathbf R }}
\newcommand{\C}{\ensuremath{ \mathbf C }}
\newcommand{\lex}{{\leq}_{lex}}
\newcommand{\slex}{{<}_{lex}}

\newcommand{\height}{\mbox{ht}}
\newcommand{\eop}{\vspace{.8cm}}
\newtheorem{theorem}{Theorem}
\newtheorem{lemma}{Lemma}
\newcommand{\bt}{\begin{theorem}}
\newcommand{\et}{\end{theorem}}
\newcommand{\bl}{\begin{lemma}}
\newcommand{\el}{\end{lemma}}
\newcommand{\pf}{{\bf Proof}.\ }
\newcommand{\bq}{\begin{eqnarray*}}
\newcommand{\eq}{\end{eqnarray*}}

\title{Polynomial growth of sumsets\\
in abelian semigroups
\footnote{
This work was begun while the authors were participants in the
Combinatorial Number Theory Workshop,
which was organized and supported by the Erd{\H o}s Center
in Budapest, Hungary.}}
\author{Melvyn B. Nathanson
\thanks{Supported in part by grants
from the PSC--CUNY Research Award Program and the
NSA Mathematical Sciences Program.}\\
Department of Mathematics\\
Lehman College (CUNY)\\
Bronx, New York 10468\\
e-mail: nathansn@alpha.lehman.cuny.edu\\
\and
Imre Z. Ruzsa\thanks{Supported in part by the Hungarian
National Foundation for Scientific Research, Grant No. T 025617.}\\
Mathematical Institute of the\\
Hungarian Academy of Sciences\\
Budapest, Pf. 127,\\
H-1364 Hungary\\
e-mail: ruzsa@math-inst.hu}
\date{}

\begin{document}
\maketitle

\begin{center}
To Michel Mend\`es-France
\end{center}

\begin{abstract}
Let $S$ be an abelian semigroup, and $A$ a finite subset of $S$.
The sumset $hA$ consists of all sums of $h$ elements of $A,$
with repetitions allowed.  Let $|hA|$ denote the cardinality of $hA$.
Elementary lattice point arguments are used to prove that
an arbitrary abelian semigroup has polynomial growth, that is,
there exists a polynomial $p(t)$ such that $|hA| = p(h)$
for all sufficiently large $h$.
Lattice point counting is also used to prove that sumsets of the form
$h_1A_1 + \cdots + h_rA_r$ have multivariate polynomial growth.
\end{abstract}

\section{Introduction}
Let $\N_0$ denote the set of nonnegative integers,
and $\N_0^k$ the set of all $k$-tuples of nonnegative integers.
Geometrically, $\N_0^k$ is the set of lattice points
in the Euclidean space $\R^k$ that lie in the nonnegative octant.

If $A$ is a finite, nonempty subset of $\N_0$,
then the {\em sumset} $hA$ is the set of all integers that can be represented
as the sum of $h$ elements of $A$, with repetitions allowed.
A classical problem in additive number theory concerns the growth of a
finite set of nonnegative integers.
For $h$ sufficiently large, the structure of the sumset $hA$
is completely determined (Nathanson~\cite{nath72f}),
and its cardinality $|hA|$ is a linear function of $h$.

If $A_1,\ldots, A_r$ are finite, nonempty subsets of $\N_0$
and if $h_1,\ldots, h_r$ are positive integers,
then $h_1A_1+\cdots + h_rA_r$ is the sumset consisting of all
integers of the form $b_1+ \cdots + b_r,$
where $b_j \in h_jA_j$ for $j = 1,\ldots,r$.
For $h_1,\ldots, h_r$ sufficiently large,
the structure of this ``linear form'' has also been completely
determined~(Han, Kirfel, and Nathanson~\cite{han-kirf-nath98}),
and its cardinality is a linear function of $h_1,\ldots, h_r$.

If $A$ is a finite, nonempty subset of $\N_0^k$,
the geometrical structure of the
sumset $hA$ is complicated, but the cardinality of $hA$
is a polynomial in $h$ of degree at most $k$ for $h$ sufficiently
large (Khovanskii~\cite{khov92}).
If the set $A$ is not contained in a hyperplane of dimension $k-1$,
then the degree of this polynomial is exactly equal to $k$.

The sets $\N_0$ and $\N_0^k$ are abelian semigroups,
that is, sets with a binary operation, called addition,
that is associative and commutative.
Let $S$ be an arbitrary abelian semigroup.
Without loss of generality, we can assume that $S$ contains
an additive identity $0$.
If $A$ is a finite, nonempty subset of $S$ and $h$ a positive integer,
we again define the sumset $hA$ as the set of all sums
of $h$ elements of $A$, with repetitions allowed.
Khovanskii~\cite{khov92,khov95} made the remarkable observation
that the cardinality of $hA$ is a polynomial
in $h$ for all sufficiently large $h$, that is,
there exists a polynomial $p(t)$ and an integer $h_0$
such that $|hA| = p(h)$ for $h \geq h_0$.
Khovanskii proved this result by constructing a finitely generated
graded module $M= \sum_{h=0}^{\infty} M_h$
over the polynomial ring $\C[t_1,\ldots,t_k]$, where $|A| = k$,
with the property that the homogeneous component $M_h$ is
a vector space over $\C$
of dimension exactly  $|hA|$ for all $h \geq 1$.
A theorem of Hilbert asserts that $\dim_{\C} M_h$ is a polynomial
in $h$ for all sufficiently large $h$, and this gives the result.

If $A_1,\ldots, A_r$ are finite, nonempty subsets of
an abelian semigroup $S$, and if $h_1,\ldots, h_r$ are positive integers,
then the ``linear form'' $h_1A_1+\cdots + h_rA_r$
is the sumset consisting of all
elements of $S$ of the form $b_1+ \cdots + b_r,$
where $b_j \in h_jA_j$ for $j = 1,\ldots,r$.
Using a generalization of Hilbert's theorem to finitely generated
modules graded by the semigroup $\N_0^r$, Nathanson~\cite{nath00c}
proved that there exists a polynomial $p(t_1,\ldots,t_r)$
such that $|h_1A_1+\cdots + h_rA_r| = p(h_1,\ldots,h_r)$
for all sufficiently large integers $h_1,\ldots, h_r$.

The purpose of this note is to give elementary combinatorial
proofs of the theorems of Khovanskii and Nathanson
that avoid the use of Hilbert polynomials.
Our arguments reduce to an easy computation about lattice points
in Euclidean space.

\section{Growth of sumsets}

We begin with some geometrical lemmas about lattice points.
Let $x = (x_1,\ldots, x_k)$ and $y = (y_1,\ldots, y_k)$
be elements of $\N_0^k$.
Define the {\em height} of $x$ by $\height(x) = \sum_{i=1}^n x_i$.
Let
\bq
\sigma(h)
& = & \{x\in \N_0^k : \height(x) = h\}  \\
& = & \{ (x_1,\ldots, x_k) \in \N_0^k : x_1+\cdots + x_k = h\}.
\eq
The set $\sigma(h)$ is a finite set of lattice points
whose cardinality is the number of ordered partitions
of $h$ as a sum of $k$ nonnegative integers, and so
\[
|\sigma(h)| = {h+k-1 \choose k-1} = \frac{h^{k-1}}{(k-1)!}
+ \frac{k h^{k-2} }{2(k-2)!}+ \cdots + 1,
\]
which is a polynomial in $h$ for fixed $k$.

We define a partial order on $\N_0^k$ by
\[
x \leq y \qquad\mbox{if $x_i \leq y_i$ for all $i = 1,\ldots, k$}.
\]
In $\N_0^2$, for example, $(2,5) \leq (4,6)$ and $(4,3) \leq (4,6)$,
but the lattice points $(2,5)$ and $(4,3)$ are incomparable.
Thus, the relation $x\leq y$ is a partial order but not a total order.
We write $x < y$ if $x \leq y$ and $x \neq y$.
If $x \leq y,$ then $x+t \leq y+t$ for all $t \in \N_0^k$.

\bl	     \label{polynomial}
Let $W$ be a finite subset of $\N_0^k$,
and let $B(h,W)$ be the set of all lattice points $x \in \sigma(h)$
such that $x \geq w$ for all $w \in W$.
Then $|B(h,W)|$ is a polynomial in $h$ for all sufficiently large $h$.
\el

\pf
Let $x = (x_1,\ldots, x_k) \in \sigma(h).$
Let $W = \{w_1,\ldots,w_m\}$,
where $w_j = (w_{1,j}, w_{2,j},\ldots, w_{k,j})\in \N_0^k$
for $j = 1,\ldots, m$.
Then $x \geq w_j$ for $j = 1,\ldots, m$ if and only if,
for all $i = 1,\ldots, k$, we have $x_i \geq w_{i,j}$  for $j=1,\ldots, m$,
that is, $x_i \geq \max\{w_{i,j}: j = 1,\ldots, m\} = w_i^*$
for $i = 1,\ldots, k$.
Define $w^* \in \N_0^k$ by $w^* = (w_1^*,\ldots, w_k^*)$.
Then
\bq
B(h,W) & = & B(h,\{w^*\}) \\
& = & \{x\in \N_0^k : \height(x) = h \mbox{ and } x \geq w^* \}  \\
& = & \{x\in \N_0^k : \height(x-w^*) = h - \height(w^*) \mbox{ and } x - w^* \geq 0 \}  \\
& = & \{y + w^*\in \N_0^k : \height(y) = h - \height(w^*) \mbox{ and } y \geq 0 \}  \\
& = & \{w^*\} + \sigma(h-\height(w^*)),
\eq
and so
\[
|B(h,W)| = |\sigma(h-\height(w^*))| = {h-\height(w^*)+k-1 \choose k-1}
\]
for $h \geq \height(w^*)$.
This completes the proof.
\eop

An {\em ideal} in an abelian semigroup is a nonempty set $I$ such that if
$x \in I$, then $x+t \in I$ for every element $t$ in the semigroup.
In the partially ordered semigroup $\N_0^k$,
a nonempty set $I$ is an ideal if and only if
$x \in I$ and $y \geq x$ imply $y \in I$.
The following result about lattice points and partial orders is
known as {\em Dickson's lemma}~\cite{cox97}.
We include a proof for completeness.

\bl				  \label{dickson}
If $I$ is a ideal in the abelian semigroup $\N_0^k$,
then there exists a finite set $W^*$ of lattice points in $\N_0^k$ such that
\[
I = \{x \in \N_0^k : x \geq w \mbox{ for some $w \in W^*$} \}.
\]
\el

\pf
The proof is by induction on the dimension $k$.
If $k=1$, then $I$ is a nonempty set of nonnegative integers, hence contains
a least integer $w$.  If $x \geq w$, then $x \in I$ since $I$ is an ideal,
and so $I = \{x\in \N_0: x \geq w\}$.

Let $k \geq 2$, and assume that the result holds for dimension $k-1$.
We shall write the lattice point $x = (x_1,\ldots,x_{k-1},x_k) \in \N_0^k$
in the form $x = (x',x_k)$, where $x' = (x_1,\ldots, x_{k-1}) \in \N_0^{k-1}$.
Define the projection map $\pi:\N_0^k\rightarrow \N_0^{k-1}$
by $\pi(x) = x'$.  Let $I' = \pi(I)$ be the image of the ideal $I$, that is,
\[
I' = \{ x' \in \N_0^{k-1} : \mbox{ $(x',x_k) \in I$ for some $x_k\in \N_0$}\}.
\]
We have $I' \neq \emptyset$ since $I \neq \emptyset$.
Let $x' \in I'$ and $y'\in \N_0^{k-1}$.
Since $x' \in I'$, there is a nonnegative integer $x_k$ such
that $(x',x_k) \in I$.
If $y' \geq x'$, then $(y',x_k) \geq (x',x_k)$ in $\N_0^k$,
and so $(y',x_k) \in I$, hence $y' \in I'$.
Thus, $I'$ is an ideal in $\N_0^{k-1}$.
Since the Lemma holds in dimension $k-1$, there is a
finite set $W' \subseteq I'$ such that $x'\in I'$
if and only if $x' \geq w'$ for some $w'\in W'$.
Associated to each lattice point $w' \in W'$ is a nonnegative integer
$x_k(w')$ such that $(w',x_k(w'))\in I$.
Let $m = \max\{x_k(w') : w'\in W'\}$ and $W_m= \{(w',m) : w' \in W'\}$.
If $w' \in W'$, then $(w',m) \geq (w',x_k(w'))$ and so $(w',m) \in I$.
Therefore, $W_m \subseteq I$.

For ${\ell}=0,1,\ldots, m-1$, we consider the set
\[
I'_{\ell} = \{x'\in \N_0^{k-1} : (x',\ell) \in I\}.
\]
If $I'_{\ell} = \emptyset$, let $W_{\ell} = \emptyset$.
If $I'_{\ell} \neq \emptyset$,
then $I'_{\ell}$ is an ideal in $\N_0^{k-1}$, and there is a finite
set $W'_{\ell}$ such that $x'\in I'_{\ell}$ if and only if $x'\geq w'$
for some $w' \in W'_{\ell}$.
Let $W_{\ell} = \{(w',\ell) : w' \in W'_{\ell}\}$.
Then $W_{\ell} \subseteq I$.
We consider the set
\[
W^* = \bigcup_{\ell =0}^{m} W_{\ell},
\]
which is a finite subset of the ideal $I$.

We shall prove that $x\in I$ if and only if $x\geq w$ for some $w \in W^*$.
If $x = (x',x_k) \in I$ and $x_k\geq m$, then $x'\in I'$, hence $x' \geq w'$
for some $w' \in W'$.   It follows that
\[
x=(x',x_k) \geq (x',m) \geq (w',m),
\]
and $(w',m) \in W_m \subseteq W^*$.

If $x = (x',\ell) \in I$ and $0 \leq  \ell < m$, then $x'\in I'_{\ell}$,
and so $x' \geq w'$ for some $w' \in W'_{\ell}$.   It follows that
\[
x = (x',\ell) \geq (w',\ell),
\]
and $(w',\ell) \in W_{\ell} \subseteq W^*$.
This completes the proof.
\eop

Let $x = (x_1,\ldots,x_k)$ and $y = (y_1,\ldots,y_k)$
be lattice points in $\N_0^k$.
We define the {\em lexicographical order} $x \lex y$
on $\N_0^k$ as follows:  $x \lex y$ if either $x = y$
or there exists $j \in \{1,2,\ldots, k\}$ such that $x_i = y_i$
for $i = 1,\ldots, j-1$ and $x_j < y_j$.
This is  a total order, so every finite, nonempty set of lattice points
contains a smallest lattice point.  For example, $(2,5) \lex (4,3) \lex (4,6)$.
If $x \lex y$, then $x+t \lex y+t$ for all $t \in \N_0^k$.
We write $x \slex y$ if $x \lex y$ and $x \neq y$

\bt  \label{khov}
Let $S$ be an abelian semigroup, and let $A$ be a finite nonempty subset of $S$.
There exists a polynomial $p(t)$ such that $|hA| = p(h)$
for all sufficiently large $h$.
\et

\pf
Let $A = \{a_1,\ldots, a_k\}$, where $|A|=k$.
We define a map $f: \N_0^k \longrightarrow S$ as follows:
If $x = (x_1,\ldots, x_k) \in \N_0^k$, then
\[
f(x) = \sum_{i=1}^k x_ia_i.
\]
This is well-defined, since each $x_i$ is a nonnegative integer
and we can add the semigroup element $a_i$ to itself $x_i$ times.
The map $f$ is a homomorphism of semigroups:  If $x,y \in \N_0^k$,
then $f(x+y) = f(x)+f(y)$.
We consider the set
\[
\sigma(h) = \{ x \in \N_0^k : \height(x) = h\}.
\]
If $x \in \sigma(h)$, then $f(x) \in hA$ and $f(\sigma(h)) = hA$.
The map $f$ is not necessarily one-to-one on the set $\sigma(h).$
For any $s \in hA$, there can be many
lattice points $x \in \sigma(h)$ such that $f(x) = s$.
However, for each $s \in hA$, there is a unique lattice point
$u_h(s) \in f^{-1}(s) \cap \sigma(h)$
that is lexicographically smallest, that is,
$u_h(s) \lex x$ for all $x \in f^{-1}(s) \cap \sigma(h)$.
Then
\[
|hA| = \left|\{u_h(s) : s \in hA\} \right|.
\]

The lattice point $x \in \N_0^k$ will be called {\em useless}
if, for $h = \height(x)$, we have $x \neq u_h(s)$ for all $s \in hA$.
Equivalently, $x \in \N_0^k$ is useless
if there exists a lattice point $u\in \sigma(\height(x))$
such that $f(u) = f(x)$ and $u \slex x.$
Let $I$ be the set of all useless lattice points in $\N_0^k$.

We shall prove that $I$ is an ideal in the semigroup $\N_0^k$.
Let $x \in I$, $\height(x) = h$, and $t \in \N_0^k$.
Since $x \in I$, there exists a lattice point $u \in \sigma(h)$
such that $f(u) = f(x)$ and $u \slex x.$
Then
\[
f(u+t) = f(u)+f(t)=f(x)+f(t) = f(x+t),
\]
\[
u+t \slex x+t,
\]
and
\[
\height(u+t) = \height(u) + \height(t) =
\height(x) + \height(t) = \height(x+t),
\]
hence
\[
u+t \in \sigma(ht(x+t)).
\]
It follows that $x+t$ is useless,
hence $x+t\in I$ and $I$ is an ideal of the semigroup $\N_0^k$.
We call $I$ the {\em useless ideal.}

By Dickson's lemma (Lemma~\ref{dickson}),
there is a finite set $W^*$ of lattice points in $\N_0^k$
such that $x \in \N_0^k$ is useless
if and only if $x \geq w$ for some $w\in W^*$.
The cardinality of the sumset $hA$ is the number of lattice points
in $\sigma(h)$ that are not in the useless ideal $I$.
For every subset $W \subseteq W^*$, we define the set
\[
B(h,W) = \{x\in \sigma(h) : x \geq w \mbox{ for all $w \in W$}\}.
\]
By the principle of inclusion-exclusion,
\[
|hA| = \sum_{W\subseteq W^*} (-1)^{|W|} |B(h,W)|.
\]
By Lemma~\ref{polynomial},
for every $W \subseteq W^*$ there is an integer $h_0(W)$
such that $|B(h,W)|$ is a polynomial in $h$ for $h \geq h_0(W)$.
Therefore, $|hA|$ is a polynomial in $h$ for all sufficiently large $h$.
This completes the proof.
\eop

\section{Growth of linear forms}

Let $k_1,\ldots,k_r$ be positive integers, and let $k = k_1 + \cdots + k_r$.
We shall write the semigroup $\N_0^k$ in the form
\[
\N_0^k = \N_0^{k_1} \times \cdots \times \N_0^{k_r},
\]
and denote the lattice point $x \in \N_0^k$ by
$x = (x_1,\ldots, x_r)$, where $x_j \in \N_0^{k_j}$ for $j = 1,\ldots, r$.
Let $h_j = \height(x_j)$ for  $j=1,\ldots,r$.
We define the {\em $r$-height} of $x$ by $\height_r(x) = (h_1,\ldots,h_r)$.
For any positive integers $h_1,\ldots, h_r$,
we consider the set
\bq
\sigma(h_1,\ldots,h_r)
& = & \{ x\in \N_0^k : \height_r(x) = (h_1,\ldots, h_r) \} \\
& = & \{ (x_1,\ldots,x_r)\in \N_0^k : \height(x_j) = h_j
\mbox{ for $j = 1,\ldots, r$}\}.
\eq
Then
\[
|\sigma(h_1,\ldots,h_r)|
= \prod_{j=1}^r |\sigma(h_j)|
= \prod_{j=1}^r {h_j+k_j-1 \choose k_j-1}
\]
is a polynomial in the $r$ variables $h_1,\ldots,h_r$ for fixed integers
$k_1,\ldots, k_r$.

\bl          \label{polynomial2}
Let $k_1,\ldots,k_r$ be positive integers, and $k = k_1 + \cdots + k_r$.
Let $W$ be a finite subset of
$\N_0^{k} = \N_0^{k_1} \times \cdots \times \N_0^{k_r}$,
and let $B(h_1,\ldots,h_r,W)$ be the set of all lattice points
$x \in \N_0^{k}$ such that $x \in \sigma(h_1,\ldots,h_r)$
and $x_j \geq w_j$
for all $w = (w_1,\ldots,w_j,\ldots, w_r) \in W$ and $j = 1,\ldots,r.$
Then $|B(h_1,\ldots,h_r,W)|$ is a polynomial in $h_1,\ldots,h_r$
for all sufficiently large integers $h_1,\ldots,h_r$.
\el

\pf
Let $x = (x_1,\ldots, x_r) \in \N_0^k$.
Let $W_j$ be the set of all lattice points $w_j \in \N_0^{k_j}$
such that there exists a lattice point $w \in W$ of the form
$w = (w_1,\ldots,w_j,\ldots,w_r)$.
Since $x \geq w$ for all $w \in W$ if and only if $x_j \geq w_j$
for all $w_j \in W_j$, it follows that the set
$B(h_1,\ldots,h_r,W)$ consists of all lattice points
$x = (x_1,\ldots,x_r) \in \N_0^k$
such that $x_j \in B(h_j,W_j)$ for all $j = 1,\ldots, r$.
Therefore,
\[
|B(h_1,\ldots,h_r,W)| = \prod_{j=1}^r |B(h_j,W_j)|.
\]
It follows from Lemma~\ref{polynomial} that $|B(h_1,\ldots,h_r,W)|$
is a polynomial in the $r$ variables $h_1,\ldots,h_r$ for all
sufficiently large integers $h_1,\ldots,h_r$.
This completes the proof.
\eop

\bt  \label{nath}
Let $S$ be an abelian semigroup, and let $A_1,\ldots, A_r$
be finite, nonempty subsets of $S$.
There exists a polynomial $p(t_1,\ldots,t_r)$ such that
$|h_1A_1 + \cdots + h_rA_r| = p(h_1,\ldots, h_r)$
for all sufficiently large integers $h_1,\ldots, h_r$.
\et

\pf
For $j = 1,\ldots, r$, let $|A_j| = k_j$
and
\[
A_j = \{a_{1,j},\ldots, a_{k_j,j}\}.
\]
Let $k = k_1 + \cdots + k_r.$
We consider lattice points
\[
x = (x_1,\ldots, x_r) \in \N_0^k = \N_0^{k_1} \times \cdots \times \N_0^{k_r},
\]
where
\[
x_j = (x_{1,j},\ldots, x_{k_j,j})\in \N_0^{k_j}.
\]
Define the semigroup homomorphism
$f: \N_0^k \rightarrow S$ as follows:
If $x = (x_1,\ldots, x_r) \in \N_0^k$,
then
\[
f(x) = \sum_{j=1}^r\sum_{i=1}^{k_j} x_{i,j}a_{i,j}.
\]
A lattice point $x \in \N_0^k$ will be called {\em $r$-useless}
if there exists a lattice point $u\in \sigma(\height_r(x))$
such that $f(u) = f(x)$ and $u \slex x.$
As in the proof of Theorem~\ref{khov},
the set $I_r$ of useless lattice points in $\N_0^k$ is an ideal.
By Lemma~\ref{dickson}, there is a finite set $W^*$
that generates $I_r$ in the sense that $x \in \N_0^k$ is $r$-useless
if and only if $x \geq w$ for some $w \in W^*$.

Let $(h_1,\ldots, h_r) \in \N_0^r$ and
\bq
\sigma(h_1,\ldots, h_r)
& = & \{(x_1,\ldots,x_r) \in \N_0^k :
\height(x_j) = h_j \mbox{ for $j=1,\ldots,r$}\}.
\eq
Then $f(\sigma(h_1,\ldots, h_r)) = h_1A_1 + \cdots + h_rA_r$,
and $|h_1A_1 + \cdots + h_rA_r|$ is the number of lattice points
in $\sigma(h_1,\ldots, h_r)$ that are not useless.
For every subset $W \subseteq W^*$, we define the set
\[
B(h_1,\ldots, h_r,W)
= \{x \in \sigma(h_1,\ldots, h_r) : x \geq w \mbox{ for all $w \in W$.}\}
\]
By the principle of inclusion-exclusion,
\[
|h_1A_1 + \cdots + h_rA_r| = \sum_{W\subseteq W^*} (-1)^{|W|} |B(h_1,\ldots, h_r,W)|.
\]
By Lemma~\ref{polynomial2},
for all sufficiently large integers $h_1,\ldots,h_r$,
the function $|B(h_1,\ldots, h_r,W)|$ is a polynomial in $h_1,\ldots,h_r$,
and so $|h_1A_1 + \cdots + h_rA_r|$ is a polynomial in $h_1,\ldots,h_r$.
This completes the proof.
\eop

{\em Remark}.
It would be interesting to describe the set of polynomials $f(t)$
such that $f(h) = |hA|$ for some finite set $A$ and sufficiently large $h$.
Similarly, one can ask for a description of the set of polynomials
$f(t_1,\ldots, t_r)$ such that
$f(h_1,\ldots,h_r) = |h_1A_1 + \cdots + h_rA_r|$, where
$A_1,\ldots, A_r$ are finite subsets of a semigroup $S$.


\begin{thebibliography}{1}

\bibitem{cox97}
D.~Cox, J.~Little, and D.~O'Shea.
\newblock {\em Ideals, Varieties, and Algorithms}.
\newblock Springer-Verlag, New York, 2nd edition, 1997.

\bibitem{han-kirf-nath98}
S.~Han, C.~Kirfel, and M.~B. Nathanson.
\newblock Linear forms in finite sets of integers.
\newblock {\em Ramanujan J.}, 2:271--281, 1998.

\bibitem{khov92}
A.~G. Khovanskii.
\newblock Newton polyhedron, {H}ilbert polynomial, and sums of finite sets.
\newblock {\em Functional. Anal. Appl.}, 26:276--281, 1992.

\bibitem{khov95}
A.~G. Khovanskii.
\newblock Sums of finite sets, orbits of commutative semigroups, and {H}ilbert
  functions.
\newblock {\em Functional. Anal. Appl.}, 29:102--112, 1995.

\bibitem{nath72f}
M.~B. Nathanson.
\newblock Sums of finite sets of integers.
\newblock {\em Amer. Math. Monthly}, 79:1010--1012, 1972.

\bibitem{nath00c}
M.~B. Nathanson.
\newblock Growth of sumsets in abelian semigroups.
\newblock {\em Semigroup Forum}, 61:149--153, 2000.

\end{thebibliography}
\end{document}